\documentclass[preprint,12pt,authoryear]{elsarticle}

\usepackage{lineno,hyperref}
\modulolinenumbers[5]

\usepackage[utf8]{inputenc}
\usepackage[utf8]{inputenc}
\usepackage[margin=0.6in]{geometry}
\usepackage{amsmath,amssymb,amsthm,graphicx,amsxtra, setspace}
\usepackage[utf8]{inputenc}
\usepackage{mathrsfs}
\usepackage{hyperref}
\usepackage{upgreek}
\usepackage{mathtools}
\usepackage{xcolor}
\usepackage{todonotes}
\allowdisplaybreaks

\usepackage{blindtext}
\newtheorem{theorem}{\noindent Theorem}[section]
\newtheorem{lemma}{\noindent Lemma}[section]
\newtheorem{corollary}{\noindent Corollary}[section]

\newtheorem{claim}{\noindent Claim}[section]
\newtheorem{remark}{\noindent Remark}[section]

\journal{Statistics \& Probability Letters}
\begin{document}

\begin{frontmatter}



\title{Study of discrete-time Hawkes process and its compensator}


\author[]{Utpal Jyoti Deba Sarma}


\author{Dharmaraja Selvamuthu}

\affiliation[]{organization={Department of Mathematics},
            addressline={Indian Institute of Technology Delhi}, 
            city={New Delhi},
            postcode={110016}, 
            state={Delhi},
            country={India}}

\begin{abstract}
The discrete-time Hawkes process (DTHP) is a sub-class of $g$-functions that serves as a discrete-time version of the continuous-time Hawkes process (CTHP). Like the CTHP, the DTHP also has the self-exciting property and its intensity depends on the entire history. In this paper, we study the asymptotic behaviour of the DTHP and its compensator. We further analyse the moment generating function (MGF) of the DTHP and obtain some bounds and convergence results on the scaled logarithmic MGF of the DTHP.
\end{abstract}




\end{frontmatter}

\section{Introduction}
The CTHP is a self-exciting point process, and was originally developed by Hawkes in 1971. \cite{hawkes1971spectra} initially studied the linear CTHP, presenting it through an immigration-birth model. The central limit theorem (CLT), law of large numbers (LLN),  large deviations principle (LDP), moderate deviations principle (MDP), etc. have been thoroughly studied by many authors (see, for instance, \cite{daley2003introduction}, \cite{bordenave2007large}, \cite{zhu2013moderate} and \cite{bacry2013some}).\par
The nonlinear CTHP, first proposed by \cite{bremaud1996stability}, are more challenging to study due to the absence of an immigration-birth representation (also called the Galton-Watson theory) and challenges in computational tractability. Nonetheless, there have been attempts to explore these processes as well.\par
CTHP are utilized across various fields, including finance, DNA modeling, seismology and neuroscience, owing to their clustering and self-exciting characteristics. In finance, these processes are particularly used for evaluating credit derivatives and modeling correlated defaults (see \cite{dassios2011dynamic} and \cite{errais2010affine}). Despite their broad applicability, CTHP may not be suitable for scenarios where data is captured in discrete time or aggregated formats. In such instances, the DTHP becomes a more relevant tool. Consequently, a discrete-time version of the CTHP has been developed to accommodate these specific data recording methods.\par
Throughout the paper, we shall define all the stochastic processes on the probability space $(\Omega,\mathcal{F},\mathbb{P})$.
\paragraph{Organisation of the paper}
Section \ref{sec 1.5} provides the definitions and literature survey of CTHP, DTHP and the compensator process. Section \ref{sec 2} focuses on deriving the intensity and compensator of the DTHP. In Section \ref{sec 3}, we present and prove the strong law of large numbers (SLLN) for DTHP. Section \ref{sec 4} discusses the SLLN, weak law of large numbers (WLLN) and CLT for the compensator of DTHP, along with its proofs. Section \ref{sec 7} examines the MGF of the DTHP and studies the convergence of the scaled logarithmic MGF. The paper concludes in Section \ref{sec 8}, summarizing the findings and discussions of the study, and outlines potential directions for future research and extensions of the work.

\section{Preliminaries and literature review}\label{sec 1.5}
In this section, we present a brief description of CTHP, DTHP, compensator and review their limit theorems. 
\paragraph{Continuous-time Hawkes process}
Consider a point process $N = \{N_t : t \geq 0\}$ on $(0,\infty)$, where $N_t$ represents the number of arrivals of events in the interval $(0, t]$. Define $\mathcal{F}_t=\sigma(N(D), D \in \mathcal{B}((0,\infty)), D \subset (0, t])$, which denotes an increasing family of $\sigma$-fields. A non-negative process $\lambda_t$ that is measurable $w.r.t$ the filtration $\mathcal{F}_t$ and satisfies
$$
\mathbb{E}\left(N(u, v] \mid \mathcal{F}_u\right) = \mathbb{E}\left(\int_u^v \lambda_s \, ds \mid \mathcal{F}_u\right),
$$
for every interval $(u, v]$, almost surely ($a.s.$), is called an intensity of $N$. A CTHP is defined as a point process $N$ that possesses an intensity of the form
$$
\lambda_t := \psi\left(\int_{0}^t g(t - s) \, dN_s\right),
$$
where $\psi: (0,\infty) \to (0,\infty)$ is a left-continuous and locally integrable function, and $g: (0,\infty) \to (0,\infty)$ is a function with a finite $L^1$ norm, i.e., $\|g\| = \int_0^{\infty} g(t) \, dt < \infty$. The process starts without prior history, denoted by $N(-\infty, 0] = 0$. The functions $\psi(\cdot)$ and $g(\cdot)$ are often called the rate function and exciting function, respectively. Depending on whether $\psi(\cdot)$ is linear or nonlinear, the CTHP is classified as linear or nonlinear respectively.
\paragraph{Discrete-time Hawkes process}\label{sec 1}
Let $\left(\beta_i\right)_{i=0}^{\infty}$ be a sequence of positive real numbers. Consider the following assumptions on this sequence.
\begin{enumerate}
  \item \label{ass1} $\sum_{i=0}^{\infty} \beta_i < 1$.
  \item \label{ass2} $\sqrt{n} \sum_{i=n}^{\infty} \beta_i \to 0$ as $n \to \infty$.
  \item \label{ass3} $\frac{1}{\sqrt{n}} \sum_{i=1}^{n} i \beta_i \to 0$ as $n \to \infty$.
  \item \label{ass4} $\sum_{i=1}^{\infty} i\beta_i < \infty$.
\end{enumerate}
Define $\{\xi_n \mid n=1,2,\ldots\}$ as a sequence of random variables having image on the set $\{0,1\}$ as follows.
\begin{itemize}
\item $\mathbb{P}(\xi_1=1)=\beta_0$; ~~
$\mathbb{P}(\xi_1=0)=1-\beta_0,$
\item
$\mathbb{P}(\xi_2=1|\xi_1)=\beta_0+\beta_1\xi_1$; ~~ 
$\mathbb{P}(\xi_2=0|\xi_1)=1-(\beta_0+\beta_1\xi_1),$
\item For $n\geq 3,$
$\mathbb{P}(\xi_n=1|\xi_1,\ldots,\xi_{n-1})=\beta_0+\sum_{i=1}^{n-1}\beta_{n-i}\xi_i$\\
$\mathbb{P}(\xi_n=0|\xi_1,\ldots,\xi_{n-1})=1-(\beta_0+\sum_{i=1}^{n-1}\beta_{n-i}\xi_i).$
\end{itemize}
Define
$$H_n=\sum_{i=1}^n\xi_i.$$\\
$\{H_n \mid n=1,2,\ldots\}$ is called the DTHP in the literature (see \cite{seol2015limit}).
\paragraph{Limit theorems on CTHP}
 In the literature, the applications of CTHP are mainly focused on the linear case. For the linear case since $\psi(\cdot)$ is linear, specifically $\psi(x) = x+ \nu$ for some $\nu > 0$ and $\|g\| < 1,$ one can examine the limit theorems using a very beautiful immigration-birth representation. \cite{daley2003introduction} proved the LLN for linear CTHP. The functional CLT for the linear multi-variate CTHP was proved under specific assumptions by \cite{bacry2013some}. This proved the CLT for linear CTHP as well. \cite{bordenave2007large} proved that if $\int_0^\infty g(t)t dt < \infty$ and $\mu \in (0,1)$, then $\left(\frac{N_t}{t}\in \cdot\right)$ satisfies the LDP, characterized by the rate function $I(x)$, where 
\begin{equation*}
I(x)= \begin{cases} x\left(\|g\|+\log \left(\frac{x}{x\|g\|+v}\right)-1\right)+\nu & \text { if } x \geq 0 \\ +\infty & \text { otherwise }\end{cases}.
\end{equation*}
The MDP for linear CTHP is studied by \cite{zhu2013moderate}.\par
 For nonlinear CTHP, \cite{bremaud1996stability} proved that under specific conditions, a stationary version of the nonlinear CTHP exists, which is unique. Furthermore, they established the convergence of a non-stationary version to equilibrium, in terms of both distribution and variation.
Afterwards, Zhu conducted extensive studies on both linear and nonlinear CTHP, as seen in his works, \cite{zhu2014process}, \cite{zhu2013central}, \cite{zhu2013moderate} and \cite{zhu2013ruin}. The CLT for nonlinear CTHP was obtained by \cite{zhu2013central}. \cite{zhu2015large} established the LDP for a particular class of nonlinear CTHP, when the exciting function takes sum of exponentials (including exponential) form. \cite{zhu2014process} also proved a level-3 LDP (also called a process-level LDP) for nonlinear CTHP with a general form of the exciting function, and subsequently derived the level-1 LDP using the contraction principle.
\paragraph{Limit theorems for DTHP}
Substantial progress has been made in studying the CTHP, however, DTHP is relatively less studied. \cite{seol2015limit} showed a WLLN for the DTHP, i.e.,
\begin{equation} \label{Eq 1}
   \frac{H_n}{n} \to \frac{\beta_0}{1-\sum_{i=1}^{\infty} \beta_i} ~ \text{as} ~ n \to \infty ~ \text{in probability}.
\end{equation}
Using the assumptions \ref{ass2} and \ref{ass3} of DTHP, \cite{seol2015limit} also derived the CLT as follows,
\begin{equation} \label{Eq 2}
\frac{H_n-n \mu}{\sqrt{n}} \to \mathcal{N}\left(0, \sigma^2\right) ~ \text{as} ~ n \to \infty ~ \text{in distribution},
\end{equation}
where
\begin{equation*}
    \mu:=\frac{\beta_0}{1-\sum_{i=1}^{\infty} \beta_i} ~ \text{and} ~ \sigma^2:=\frac{(1-\mu)\mu}{\left(1-\sum_{j=1}^{\infty} \beta_j\right)^2}.
\end{equation*}
\cite{wang2022limit} studied a marked DTHP introduced in \cite{xu2021self}, focusing on a more general exciting function and deduced the
LLN and CLT for the model. The MDP and LDP for the model are studied in \cite{wang2023large}.
\paragraph{Other related literature on CTHP}
Authors have recently started studying different modifications of CTHP. Once such modification is the extended inverse CTHP. The inverse CTHP is characterized by a fixed intensity (in contrast to the CTHP that has stochastic intensity), and jump size, which is stochastic in nature. An extended inverse CTHP is formed through the combination of a CTHP and an inverse CTHP. \cite{selvamuthu2023limit} analyzed the limiting properties of this extended inverse CTHP, with a general exciting function, and they derived the CLT the LLN.\par
\cite{seol2019limit} examined the inverse Markovian CTHP, combining features of various self-exciting processes. He specifically conducted an asymptotic analysis for a variant of the inverse Markovian CTHP.
\paragraph{Compensator of a point process}
For a point process $N$, its compensator, denoted by $\Lambda_t$,  is defined as the unique increasing $\mathcal{F}_t$-predictable process starting from zero, and satisfying the equation 
$$N_t=M_t+\Lambda_t,$$ 
$a.s.$ for all $t \geq 0$, where $M_t$ is a martingale. The Doob's decomposition theorem assures the presence of such an $M_t$.\par
In the literature, authors have studied the compensator of point processes, especially within the context of CTHP. The compensator is crucial in analyzing the dynamical behaviors of point processes because the point process and its compensator is separated only by a martingale. 
For the importance of compensator, refer to \cite{daley2003introduction}, \cite{guo2008intensity} and \cite{ivanoff2007compensator}.
\cite{SEOL2017165} examined the compensator of linear CTHP, and proved the LLN, CLT and the LDP for the compensator process.
\section{Intensity and compensator of DTHP}\label{sec 2}
In this section, we introduce the intensity of the DTHP. Then, we derive the compensator of the DTHP.
\paragraph{Intensity}
Let $\mathcal{F}_n=\sigma(\xi_1,\ldots,\xi_n)$, i.e., $\mathcal{F}_n$ is the $\sigma$-field generated by $\xi_1,\ldots,\xi_n$.\\
From the DTHP described in Section \ref{sec 1}, we can conclude that $\lambda_n:=\beta_0 + \sum_{i=1}^{n-1} \beta_{n-i}\xi_i$ is the intensity of the DTHP $\{H_n \mid n=1,2,\ldots\}$. The sequence $\left(\beta_i\right)_{i=0}^\infty$ plays the role of the exciting function for the DTHP. The arrival of each event at discrete time $n$ increases the intensity, establishing the self-exciting nature of the DTHP. Also, the intensity of the DTHP depends on its entire history. All these characteristics make it a discrete-time equivalent of the CTHP.
\paragraph{Compensator}
Now, we show that the DTHP $\{H_n \mid n=1,2,\ldots\}$ defined in Section \ref{sec 1} is a sub-martingale $w.r.t.$ the filtration $\{\mathcal{F}_n \mid n=1,2,\ldots\}$. Then we use Doob's decomposition theorem to get the compensator of $\{H_n \mid n=1,2,\ldots\}$.
\begin{lemma}
 The DTHP is a sub-martingale $w.r.t.$ the filtration $\{\mathcal{F}_n \mid n=1,2,\ldots\}$.
\end{lemma}

\begin{proof}
From the definition of $\{H_n \mid n=1,2,\ldots\}$,
$$\left | H_n \right | = \left|\sum_{i=1}^n \xi_i \right| \leq n.$$
Thus, $H_n$ is bounded for each $n=1,2,\ldots$, which implies $H_n$ is integrable for each $n=1,2,\ldots.$ Also, from the definition, $H_n$ is $\mathcal{F}_n$-measurable.\\
Finally
\begin{align*}
    \mathbb{E}(H_{n+1}\mid\mathcal{F}_n)= \mathbb{E}(H_n+\xi_{n+1}\mid\mathcal{F}_n)
    =H_n+\mathbb{E}(\xi_{n+1}\mid \mathcal{F}_n)
    \geq H_n.
    \end{align*}
Hence, the DTHP is a sub-martingale $w.r.t.$ the filtration $\{\mathcal{F}_n \mid n=1,2,\ldots\}$.
\end{proof}
Since the DTHP is a sub-martingale, by Doob's decomposition, it can be written as a sum of an increasing predictable process and a martingale, uniquely, i.e.,
\begin{equation}\label{eq 1}
    H_n=M_n+\Lambda_n,
    \end{equation}
    where $M_n=\sum_{i=1}^n \left[\xi_i-\mathbb E(\xi_i \mid \mathcal{F}_{i-1}) \right]$ is a martingale and $\Lambda_n=\sum_{i=1}^n \mathbb E(\xi_i \mid \mathcal{F}_{i-1})$ is an increasing predictable process.
The fact that $M_n$ is a martingale and $\Lambda_n$ is an increasing predictable process is verified as follows.
\begin{claim}\label{claim 1}
    $\{M_n \mid n=1,2,\ldots\}$ is a martingale $w.r.t.$ the filtration $\{\mathcal{F}_n \mid n=1,2,\ldots\}$.
\end{claim}
\begin{proof}
    From the definition of $M_n$, $M_n$ is $\mathcal{F}_n$-measurable.\\
    Also
    \begin{align*}
        \left| M_n \right| = \left| \sum_{i=1}^n \left(\xi_i-\mathbb E(\xi_i \mid \mathcal{F}_{i-1})\right) \right|\leq \sum_{i=1}^n \left| \xi_i - \mathbb E(\xi_i \mid\mathcal{F}_{i-1})\right|\leq n.
    \end{align*}
    Hence, $M_n$ is bounded for each $n=1,2,\ldots$, which implies that $M_n$ is integrable for each $n=1,2,\ldots.$\\
    Finally
    \begin{align*}
        \mathbb{E}(M_{n+1}\mid \mathcal{F}_n)=\mathbb{E}\left(\sum_{i=1}^{n+1} \left[\xi_i-\mathbb E(\xi_i \mid \mathcal{F}_{i-1})\right] \mid \mathcal{F}_n\right)
        = \sum_{i=1}^n \left[\xi_i-\mathbb E(\xi_i \mid \mathcal{F}_{i-1})\right] + \mathbb{E}\left(\xi_{n+1}-\mathbb{E}(\xi_{n+1}\mid \mathcal{F}_n)\mid \mathcal{F}_n\right)= M_n.
    \end{align*}
    Hence, $\{M_n \mid n=1,2,\ldots\}$ is a martingale $w.r.t.$ the filtration $\{\mathcal{F}_n \mid n=1,2,\ldots\}$.
\end{proof}
\begin{remark}\label{rem}
    We assume for convenience that $M_0=0$.
\end{remark}
\begin{claim}
    $\{\Lambda_n \mid n=1,2,\ldots\}$ is an increasing predictable process.
\end{claim}
\begin{proof}
    From the definition of $\Lambda_n$, we have
    \begin{align*}
        \Lambda_{n+1}-\Lambda_n&=\sum_{i=1}^{n+1} \mathbb E(\xi_i \mid \mathcal{F}_{i-1}) - \sum_{i=1}^{n} \mathbb E(\xi_i \mid \mathcal{F}_{i-1})
        =\mathbb E(\xi_{n+1} \mid \mathcal{F}_n)
        \geq 0.
\end{align*}
    Also, from the definition, $\Lambda_n$ is $\mathcal{F}_{n-1}$-measurable, which implies it is a predictable process.\\
    Hence, $\{\Lambda_n \mid n=1,2,\ldots\}$ is an increasing predictable process. 
\end{proof}
This process $\{\Lambda_n \mid n=1,2,\ldots\}$ is the compensator of the DTHP $\{H_n \mid n=1,2,\ldots\}$.
\section{SLLN for DTHP}\label{sec 3}
In this section, we obtain the SLLN for the DTHP $\{H_n \mid n=1,2,\ldots\}$ defined in Section \ref{sec 1}.
\begin{theorem}\label{thm 1}
 Under assumptions \ref{ass1}, \ref{ass2}, \ref{ass3} and \ref{ass4},
$$\frac{H_n}{n} \to \frac{\beta_0}{1 - \sum_{i=1}^{\infty} \beta_i}, ~ \text{as} ~ n \to \infty ~ a.s..$$
\end{theorem}
\paragraph{Proof of Theorem \ref{thm 1}}
    From the definition of $\{\xi_n \mid n=1,2,\ldots\}$,
\begin{align*}
    \mathbb{E}(\xi_i \mid \mathcal{F}_{i-1})&=0\cdot\mathbb{P}(\xi_i=0\mid\mathcal{F}_{i-1})+1\cdot\mathbb{P}(\xi_i=1\mid\mathcal{F}_{i-1})
    =\beta_0 + \sum_{i=1}^{n-1} \beta_{n-i}\xi_i.
\end{align*}
Now
\begin{align*}
    M_n = \sum_{i=1}^n \left[\xi_i - \mathbb{E}(\xi_i | \mathcal{F}_{i-1}) \right]
    = H_n  - \beta_0 n - \sum_{i=1}^n \sum_{j=1}^{i-1} \beta_{i-j}\xi_j.
    \end{align*}
    After some computation, the expression $\sum_{i=1}^n \sum_{j=1}^{i-1} \beta_{i-j}\xi_j$ can be rewritten as $\sum_{i=1}^{n-1} \sum_{j=1}^{n-i} \beta_j \xi_i.$\\
    Therefore,
    \begin{align*}
    M_n= H_n - \beta_0 n - \sum_{i=1}^{n-1} \sum_{j=1}^{n-i} \beta_j \xi_i
    = \left(1- \sum_{j=1}^{\infty}\beta_j\right)H_n-\beta_0n+\left(\sum_{j=1}^{\infty}\beta_j\right)H_n - \sum_{i=1}^{n-1}\sum_{j=1}^{n-i}\beta_j\xi_i.
    \end{align*}
    It can be rewritten as
    \begin{equation}\label{eq 2}
    M_n =  \left(1- \sum_{j=1}^{\infty}\beta_j\right)H_n-\beta_{0}n+ \zeta_n,
\end{equation}
where $$\zeta_n=\left(\sum_{j=1}^{\infty}\beta_j\right)H_n - \sum_{i=1}^{n-1}\sum_{j=1}^{n-i}\beta_j\xi_i.$$\\
\begin{lemma} \label{lemma 4.1}
    \begin{equation*}
 \frac{\zeta_n}{n}\to 0, ~ \text{as} ~ n\to\infty ~ a.s..
\end{equation*}
\end{lemma}
\begin{proof}
\begin{align*}
    \zeta_n&=\left(\sum_{j=1}^{\infty}\beta_j\right)H_n - \sum_{i=1}^{n-1}\sum_{j=1}^{n-i}\beta_j\xi_i\\
    &=\left(\sum_{j=1}^{\infty}\beta_j\right)H_n - \left[\left(\sum_{i=1}^{n-1}\beta_i\right)\xi_1+\left(\sum_{i=1}^{n-2}\beta_i\right)\xi_2+\ldots+\left(\sum_{i=1}^{2}\beta_i\right)\xi_{n-2}+ \beta_1\xi_{n-1}\right]\\
    &= \beta_1\xi_n+\beta_2\left(\sum_{i=0}^{1}\xi_{n-i}\right)+\beta_3\left(\sum_{i=0}^2\xi_{n-i}\right)+\ldots+\beta_n\left(\sum_{i=0}^{n-1}\xi_{n-i}\right)+\left(\sum_{j=n}^{\infty}\beta_j\right)H_n\\
    &\leq \beta_1+2\beta_2+3\beta_3+\ldots+n\beta_n+n\left(\sum_{j=n+1}^{\infty}\beta_j\right) 
    \leq \sum_{j=1}^{\infty}j\beta_j.
\end{align*}
Hence, we have
   \begin{equation}\label{ineq zeta}
   \frac{\zeta_n}{n} \leq  \frac{\sum_{j=1}^{\infty}j\beta_j}{n} ~ \forall ~ n=1,2,\ldots.
   \end{equation}
   By assumption \ref{ass4} of the DTHP described in Section \ref{sec 1}, and the fact that $\zeta_n\geq0 ~ \forall ~ n=1,2,\ldots$, we get   \begin{equation*}
 \frac{\zeta_n}{n}\to 0, ~ \text{as} ~ n\to\infty ~ a.s..
\end{equation*}
\end{proof}
\begin{lemma}\label{lemma 4.2}
    \begin{equation*}
    \frac{M_n}{n} \to 0, ~ \text{as} ~ n\to\infty ~ a.s..
    \end{equation*}
\end{lemma}
\begin{proof}
Let $$A_n=M_n-M_{n-1}=\xi_{n}-\mathbb{E}(\xi_{n}\mid \mathcal{F}_{n-1}).$$
This implies 
$$|A_n|\leq1 ~ \text{and} ~ \mathbb{E}(A_n \mid \mathcal{F}_{n-1})=0.$$
Also, for any $n_1<n_2<\ldots < n_k \in \mathbb{N},$
    \begin{align*}
        \mathbb{E}(A_{n_1}A_{n_2}\ldots A_{n_k})
        =\mathbb{E}[\mathbb{E}(A_{n_1}A_{n_2}\ldots A_{n_k})\mid \mathcal{F}_{n_k-1}]=\mathbb{E}[A_{n_1}A_{n_2}\ldots A_{n_{k-1}}\mathbb{E}(A_{n_k}\mid \mathcal{F}_{n_k-1})]
        =0.
    \end{align*}
    Let 
    $$B_n=\sum_{i=1}^n \frac{A_i}{i}.$$
    By definition, $B_n$ is integrable for each $n=1,2,\ldots$, and $B_n$ is $\mathcal{F}_n$-measurable.\\
    Also,
    \begin{align*}
        \mathbb{E}(B_{n+1}\mid \mathcal{F}_n)=\sum_{i=1}^n \frac{A_i}{i}+\mathbb{E}\left(\frac{A_{n+1}}{n+1} \Bigm| \mathcal{F}_n\right)= B_n.
    \end{align*}
    Thus, $\{B_n \mid n=1,2,\ldots\}$ is a martingale $w.r.t$ the filtration $\{\mathcal{F}_n \mid n=1,2,\ldots\}$. 
    Now,
    \begin{align*}
        \mathbb{E}(B_{n}^2)=\sum_{i=1}^n\frac{\mathbb{E}(A_{i}^2)}{i^2}+2\sum_{i<j}\frac{\mathbb{E}(A_{i}A_{j})}{ij}\leq \sum_{i=1}^n \frac{1}{i^2}.
    \end{align*}
    Thus, $$ \sup_{n\in \mathbb{N}} \mathbb{E}(B_{n}^2) < +\infty.$$
    Using Jensen's inequality, $\{B_{n}^2 \mid n=1,2,\ldots\}$ is a sub-martingale $w.r.t$ the filtration $\{\mathcal{F}_n \mid n=1,2,\ldots\}.$ Hence, by Doob's Martingale Convergence Theorem, there exists an integrable random variable $B$ such that 
\begin{equation*}
    B_{n}^2 \to B ~ \text{as} ~ n \to \infty ~ a.s..
\end{equation*}
This implies
\begin{equation*}
    |B_{n}| \to \sqrt{B} ~ \text{as} ~ n \to \infty ~ a.s..
\end{equation*}
Thus,
\begin{equation*}
    \Big|\sum_{i=1}^{\infty}\frac{A_i}{i}\Big|<\infty ~ a.s., ~ \text{which implies} ~ \sum_{i=1}^{\infty}\frac{A_i}{i}<\infty ~ a.s..
\end{equation*}
Hence, by Kronecker's lemma,
\begin{equation*}
    \frac{1}{n}\sum_{i=1}^n i \frac{A_i}{i} \to 0 ~ \text{as} ~ n \to \infty ~ a.s.,
\end{equation*}
which implies, using Remark \ref{rem}
\begin{equation}\label{con 4}
    \frac{M_n}{n} \to 0 ~ \text{as} ~n \to \infty ~ a.s..
\end{equation}
\end{proof}

\begin{corollary}
    Since convergence $a.s.$ implies convergence in probability, hence, from convergence relation (\ref{con 4}) we can conclude that
    \begin{equation}\label{con 5}
    \frac{M_n}{n} \to 0, ~ \text{as} ~ n \to \infty ~ \text{in probability}.
\end{equation}
\end{corollary}\par
From Equation (\ref{eq 2}), we get
    \begin{equation*}
    \frac{M_n}{n}=\left(1- \sum_{j=1}^{\infty}\beta_j\right)\frac{H_n}{n}-\beta_0+ \frac{\zeta_n}{n}.
\end{equation*}
Using Lemma \ref{lemma 4.1} and \ref{lemma 4.2}, letting $n\to\infty$ on both sides we get
\begin{equation}\label{eq 4.4}
    \frac{H_{n}}{n}\to \frac{\beta_0}{1- \sum_{j=1}^{\infty}\beta_j}, ~ \text{as} ~ n\to\infty ~ a.s..
    \end{equation}
\section{LLN and CLT for the compensator of DTHP}\label{sec 4}
In this section, we obtain the SLLN, WLLN and CLT for the compensator of the DTHP.
\begin{theorem}[SLLN]\label{thm 2}
    $$\frac{\Lambda_n}{n} \to \frac{\beta_0}{1 - \sum_{i=1}^{\infty} \beta_i}, ~ \text{as} ~ n\to\infty ~ a.s..$$
\end{theorem}
\begin{theorem}[WLLN]\label{thm 3}
    \begin{equation*}
    \frac{\Lambda_n}{n} \to \frac{\beta_0}{1 - \sum_{i=1}^{\infty} \beta_i}, ~ \text{as} ~ n \to \infty ~ \text{in probability}.
\end{equation*}
\end{theorem}
\begin{theorem}[CLT]\label{thm 4}
    \begin{equation*}
        \frac{\Lambda_n - \mu n}{\sqrt{n}} \to \mathcal{N}\left(0, \left(\sum_{i=1}^{\infty} \beta_i\right)^2 \sigma^2 \right), ~ \text{as} ~n \to \infty ~ \text{in distribution}.
\end{equation*}
where
$$
\mu := \frac{\beta_0}{1 - \sum_{i=1}^{\infty} \beta_i} ~ \text{and} ~ \sigma^2 := \frac{\mu(1-\mu)}{\left(1-\sum_{i=1}^{\infty} \beta_i\right)^2}.
$$
\end{theorem}
\paragraph{Proof of Theorem \ref{thm 2}}
From Equation (\ref{eq 1}), we have
\begin{equation}\label{Eq lambda}
    \frac{\Lambda_n}{n}=\frac{H_n}{n}-\frac{M_n}{n}.
\end{equation}
Hence, using convergence relations (\ref{con 4}) and (\ref{eq 4.4}), we conclude that
$$\frac{\Lambda_n}{n} \to \frac{\beta_0}{1 - \sum_{i=1}^{\infty} \beta_i},  ~ \text{as} ~ n \to \infty ~ a.s..$$
\paragraph{Proof of Theorem \ref{thm 3}}
From Equation (\ref{Eq lambda}), and using WLLN for $\{H_n \mid n=1,2,\ldots\}$, i.e., convergence relation (\ref{Eq 1}) and convergence relation (\ref{con 5}), we conclude that
$$\frac{\Lambda_n}{n} \to \frac{\beta_0}{1 - \sum_{i=1}^{\infty} \beta_i}, ~ \text{as} ~ n\to\infty ~ \text{in probability}.$$

\paragraph{Proof of Theorem \ref{thm 4}}
From the definition of compensator of DTHP,
\begin{align*}
    \Lambda_n =\sum_{i=1}^n \mathbb E(\xi_i \mid \mathcal{F}_{i-1})
    =\sum_{i=1}^n \left(\beta_0+\sum_{j=1}^{i-1} \beta_{i-j}\xi_j\right)
    =n\beta_0+\sum_{i=1}^n\sum_{j=1}^{i-1} \beta_{i-j}\xi_j
    =n\beta_0+\sum_{i=1}^{n-1}\sum_{j=1}^{n-i} \beta_{j}\xi_i.
\end{align*}
Using Equations (\ref{eq 1}) and (\ref{eq 2}), we get
$$\Lambda_n =  \left(\sum_{j=1}^{\infty}\beta_j\right)H_n+\beta_0n- \zeta_n.$$\\
From Inequality (\ref{ineq zeta}) in the proof of Lemma \ref{lemma 4.1}, we have the following
   \begin{equation}\label{ineq 4.8}
   \frac{\mathbb{E}(\zeta_n)}{\sqrt{n}}\leq \frac{\sum_{j=1}^{\infty}j\beta_j}{\sqrt{n}}.
   \end{equation}
Using the fact that $\zeta_n\geq 0$ and $\sum_{j=1}^{\infty}j\beta_j<\infty,$  we get
$$\frac{\mathbb{E}(\zeta_n)}{\sqrt{n}} \to 0, ~ \text{as} ~ n\to\infty ~ a.s..$$
Consider the following expression,
\begin{align*}
    \frac{\Lambda_n -n\frac{\beta_0}{1 - \sum_{i=1}^{\infty} \beta_i}}{\sqrt{n}}
    &=\frac{\left(\sum_{j=1}^{\infty}\beta_j\right)H_n+\beta_0n- \zeta_n - n\frac{\beta_0}{1 - \sum_{j=1}^{\infty} \beta_i}}{\sqrt{n}}\\
    &=\frac{\left(\sum_{j=1}^{\infty}\beta_j\right)H_n-n\frac{\left(\sum_{j=1}^{\infty} \beta_j\right)\beta_0}{1 - \sum_{j=1}^{\infty} \beta_i}-\zeta_n}{\sqrt{n}}\\
    &=\frac{\left(\sum_{j=1}^{\infty}\beta_j\right)\left(H_n-n\frac{\beta_0}{1 - \sum_{j=1}^{\infty} \beta_j}\right)-\zeta_n}{\sqrt{n}}.
\end{align*}
Hence, we get
\begin{equation*}
\frac{\Lambda_n -n\frac{\beta_0}{1 - \sum_{i=1}^{\infty} \beta_i}}{\sqrt{n}}=\frac{\left(\sum_{j=1}^{\infty}\beta_j\right)\left(H_n-n\frac{\beta_0}{1 - \sum_{j=1}^{\infty} \beta_j}\right)}{\sqrt{n}}-\frac{\zeta_n}{\sqrt{n}}.
\end{equation*}
Using the Slutsky's theorem, CLT of $\{H_n \mid n=1,2,\ldots\}$, i.e., convergence relation (\ref{Eq 2}), and the fact that $\frac{\mathbb{E}(\zeta_n)}{\sqrt{n}} \to 0$ as $n\to\infty ~ a.s.$, we get
$$
\frac{\Lambda_n-n\mu}{\sqrt{n}} \to \mathcal{N}\left(0, \frac{(\sum_{j=1}^{\infty} \beta_j)^2\mu(1-\mu)}{\left(1-\sum_{j=1}^{\infty} \beta_j\right)^2}\right), ~ \text{as} ~ n \to \infty ~ \text{in distribution}.
$$
\section{Scaled logarithmic moment generating function of DTHP}\label{sec 7}
In this section, we derive some bounds of the scaled logarithmic MGF of the DTHP, and also prove some important convergence results.\\
Consider the MGF of $H_n$.
\begin{align*}
\mathbb{E}(e^{tH_n}) &= \mathbb{E}\left( e^{t(H_{n-1} + \xi_n)} \right) = \mathbb{E}\left( \mathbb{E}(e^{tH_{n-1} +{t\xi_n}} \mid \mathcal{F}_{n-1} \right) = \mathbb{E}\left(e^{tH_{n-1}}\mathbb{E}(e^{t\xi_n}\mid \mathcal{F}_{n-1})\right)\\
&= \mathbb{E}\left(e^{tH_{n-1}}\left[\mathbb{P}(\xi_n=0\mid \mathcal{F}_{n-1})+e^t\mathbb{P}(\xi_n=1\mid \mathcal{F}_{n-1})\right]\right)\\
&= \mathbb{E}\left(e^{tH_{n-1}}-e^{tH_{n-1}}\mathbb{P}(\xi_n=1\mid \mathcal{F}_{n-1})+e^{tH_{n-1}}e^t\mathbb{P}(\xi_n=1\mid \mathcal{F}_{n-1})\right)\\
&= \mathbb{E}\left(e^{tH_{n-1}}+\mathbb{P}(\xi_n=1\mid \mathcal{F}_{n-1})e^{tH_{n-1}}(e^t-1)\right)\\
&= \mathbb{E}(e^{tH_{n-1}})+(e^t-1)\mathbb{E}(e^{tH_{n-1}}\mathbb{P}(\xi_n=1\mid \mathcal{F}_{n-1})).
\end{align*}
Hence, we obtain the relation
\begin{align}\label{eq 7}
\mathbb{E}(e^{tH_{n}})=\mathbb{E}(e^{tH_{n-1}})+(e^t-1)\mathbb{E}(e^{tH_{n-1}}\mathbb{P}(\xi_n=1\mid \mathcal{F}_{n-1})).
\end{align}
Define the scaled logarithmic MGF of the DTHP as \begin{equation}\label{eq gamma} 
\Gamma_n(t)=\frac{1}{n}\log \mathbb{E}(e^{tH_n}).
\end{equation}
We now obtain the bounds for $\Gamma_n$.
\begin{theorem}\label{thm 5}
    For each $t \in \mathbb R$, $\Gamma_n(t)$ is pointwise bounded, i.e., there exists two functions $\psi(t): \mathbb R \to \mathbb R$ and $\phi(t): \mathbb R \to \mathbb R$ such that 
$$\phi(t) \leq \Gamma_n(t)\leq \psi(t).$$
\end{theorem}
\paragraph{Proof of Theorem \ref{thm 5}}
We divide the proof into four lemmas as follows.
\begin{lemma}
    For $t>0$, $\Gamma_n(t)$ is bounded above.
\end{lemma}
\begin{proof}
From Equation \eqref{eq 7}, for $t>0$
\begin{align*}
\mathbb{E}(e^{tH_{n}})&=\mathbb{E}(e^{tH_{n-1}})+(e^t-1)\mathbb{E}\left(e^{tH_{n-1}}(\beta_0+\sum_{i=1}^{n-1}\beta_{n-i}\xi_i)\right)\\
&\leq \mathbb{E}(e^{tH_{n-1}})+(e^t-1)\mathbb{E}\left(e^{tH_{n-1}}(\beta_0+ \ldots +\beta_{n-1})\right)\\
&= \mathbb{E}(e^{tH_{n-1}})+(e^t-1)(\beta_0+ \ldots +\beta_{n-1})\mathbb{E}\left(e^{tH_{n-1}}\right)\\
&= \mathbb{E}(e^{tH_{n-1}})\left[1+(e^t-1)(\beta_0+ \ldots +\beta_{n-1})\right]\\
&= \mathbb{E}(e^{tH_{n-1}})\left[e^t(\beta_0+ \ldots +\beta_{n-1})+(1-(\beta_0+ \ldots +\beta_{n-1}))\right]\\
&\leq e^{t} \mathbb{E}(e^{tH_{n-1}}).
\end{align*}
The last inequality follows from the following.
Suppose 
$g(t)=e^t-[e^t(\beta_0 + \ldots + \beta _{n-1})+[1-(\beta_0 + \ldots + \beta _{n-1})].$
Then, for $t>0,
    g(t)=(e^t-1)\left(1-(\beta_0 + \ldots + \beta _{n-1})\right) \geq 0.$\par
Continuing recursively,
\begin{align*}
\mathbb{E}(e^{tH_{n}})\leq e^{t} \mathbb{E}(e^{tH_{n-1}})
\leq e^{2t} \mathbb{E}(e^{tH_{n-2}})\leq \ldots \leq e^{(n-1)t} \mathbb{E}(e^{tH_1})= e^{(n-1)t}((1-\beta_0)+e^t\beta_0))\leq e^{nt}.
\end{align*}
Thus, for $t>0$
\begin{align*}
     \Gamma_n(t)=\frac{\log\mathbb{E}(e^{tH_n})}{n} \leq {\frac{\log e^{nt}}{n}}=\frac{nt}{n}=t.
    \end{align*}
\end{proof}
\begin{lemma}
    For $t>0$, $\Gamma_n(t)$ is bounded below.
\end{lemma}
\begin{proof}
For $t>0$, from Equation (\ref{eq 7}) we have,
\begin{align*}
\mathbb{E}(e^{tH_{n}})&=\mathbb{E}(e^{tH_{n-1}})+(e^t-1)\mathbb{E}\left(e^{tH_{n-1}}(\beta_0+\sum_{i=1}^{n-1}\beta_{n-i}\xi_i)\right)\\
&\geq \mathbb{E}(e^{tH_{n-1}})+(e^t-1)\mathbb{E}(e^{tH_{n-1}}\beta_0)\\
&= \mathbb{E}(e^{tH_{n-1}})\left(1+\beta_0(e^t-1)\right)=\mathbb{E}(e^{tH_{n-1}})\mathbb{E}(e^{tH_{1}}).
\end{align*}
Recursively,
\begin{align*}
\mathbb{E}(e^{tH_{n}}) &\geq \mathbb{E}(e^{tH_{n-1}})\mathbb{E}(e^{tH_{1}})\geq \ldots \geq \mathbb{E}(e^{tH_{1}})^{n-1}\mathbb{E}(e^{tH_1})= (1-\beta_0+\beta_0e^t)^{n}.
\end{align*}
Hence, for $t>0$
\begin{align*}
\Gamma_n(t)=\frac{\log \mathbb{E}(e^{tH_n})}{n}\geq \frac{\log(1-\beta_0+\beta_0e^t)^n}{n}= \frac{n\log(1-\beta_0+\beta_0e^t)}{n}=\log(1-\beta_0+\beta_0e^t).
\end{align*}
\end{proof}
\begin{lemma}
    For $t<0$, $\Gamma_n(t)$ is bounded above.
\end{lemma}
\begin{proof}
From equation \eqref{eq 7}, for $t<0$
\begin{align*}
    \mathbb{E}(e^{tH_{n}})=\mathbb{E}(e^{tH_{n-1}})+(e^t-1)\mathbb{E}\left(e^{tH_{n-1}}(\beta_0+\sum_{i=1}^{n-1}\beta_{n-i}\xi_i)\right) \leq \mathbb{E}(e^{tH_{n-1}}).
\end{align*}
Continuing recursively,
\begin{align*}
    \mathbb{E}(e^{tH_{n}}) &\leq \mathbb{E}(e^{tH_{n-1}})\leq \ldots \leq\mathbb{E}(e^{tH_1})= (1-\beta_0)+\beta_0 e^t.
\end{align*}
Hence, for $t<0$,
\begin{align*}
\Gamma_n(t)=\frac{\log\mathbb{E}(e^{tH_{n}})}{n}\leq \frac{\log((1-\beta_0)+\beta_0 e^t )}{n}\leq 0.
\end{align*}
\end{proof}
\begin{lemma}
    For $t<0$, $\Gamma_n(t)$ is bounded below.
\end{lemma}
\begin{proof}
For any $n=1,2,\ldots,$ the MGF of $H_n$ can also be written as
\begin{align*}
    \mathbb{E}(e^{tH_{n}})=c_{0,n}+c_{1,n} e^t +\ldots + c_{n,n} e^{nt}.
\end{align*}
where
$$c_{r,n}=\sum_{\substack{a_1+\ldots +a_n=r\\ a_i=0 ~ \text{or} ~ 1}} \mathbb{P}(\xi_1=a_1, \xi_2=a_2, \ldots, \xi_n=a_n),~ \sum_{r=0}^n c_{r,n}=1, ~ 0\leq r \leq n.$$
Using multiplicative rule of probability, we get
\begin{eqnarray*}
c_{0,n}=\mathbb{P}(\xi_1=0, \xi_2=0, \ldots, \xi_n=0)
&=&\mathbb{P}(\xi_1=0)\cdot\mathbb{P}(\xi_2=0\mid \xi_1=0)\cdot \\
& & \ldots \cdot \mathbb{P}(\xi_n=0\mid \xi_1=0,\xi_2=0,\ldots,\xi_{n-1}=0)\\
 &= & (1-\beta_0)^n.
\end{eqnarray*}
Similarly,
\begin{align*}
c_{n,n}=\beta_0(\beta_0+\beta_1)\ldots(\beta_0+\ldots+\beta_n).
\end{align*}
For $t<0$,
\begin{align*}
&\mathbb{E}(e^{tH_{n}})=c_{0,n}+c_{1,n} e^t +\ldots + c_{n,n} e^{nt} \geq c_{0,n},
\end{align*}
Hence, for $t<0$
\begin{align*}
\Gamma_n(t) = \frac{\log\mathbb{E}(e^{tH_{n}})}{n}\geq \frac{\log c_{0,n}}{n}= \frac{\log(1-\beta_0)^n}{n}= \log(1-\beta_0).
\end{align*}
\end{proof}
\subsection{Convergence of \texorpdfstring{$\Gamma_n(t)$}{Gamma\_n(t)}}

In this section, we prove that  for $t<0$, $\Gamma_n(t)$ converges to some function of $t$. First, we will prove the following lemma.
\begin{lemma}\label{lemma 10}
    For $t<0$, $\Gamma_n(t)$ is monotonically decreasing.
\end{lemma}
\begin{proof}
From equation \eqref{eq 7}, for any $n=2,3,\ldots$ and $t<0$
\begin{align*}
 &\mathbb{E}(e^{tH_{n}}) < \mathbb{E}(e^{tH_{n-1}}).
 \end{align*}
 Hence, for $t<0$,
 \begin{align*}
 \Gamma_n(t)=\frac{\log\mathbb{E}(e^{tH_{n}})}{n}<\frac{\log\mathbb{E}(e^{tH_{n-1}})}{n}<\frac{\log\mathbb{E}(e^{tH_{n-1}})}{n-1}=\Gamma_{n-1}(t).
 \end{align*}
\end{proof}
By Theorem \ref{thm 5}, $\Gamma_n(t)$ is bounded for all $t\in\mathbb{R}$, and by Lemma \ref{lemma 10}, $\Gamma_n(t)$ is monotonically decreasing for  $t<0$, hence from monotone convergence theorem, we conclude that for $t<0$, $\Gamma_n(t)$ converges, i.e., for $t<0$
\begin{equation*}
    \lim_{n \to \infty} \Gamma_n(t)=\Gamma(t) .
\end{equation*}
\section{Conclusion and future work} \label{sec 8}
In this paper, we have extended the work of \cite{seol2015limit} to prove the SLLN for DTHP $\{H_n \mid n=1,2,\ldots \}$. We have identified the compensator of the DTHP and analyzed its asymptotic behaviour. We have shown the WLLN, the SLLN and the CLT for the compensator of the DTHP. We have further analyzed the MGF and scaled logarithmic MGF of the DTHP obtaining tight bounds and convergence results. The analysis on the scaled logarithmic MGF is important because these results are necessary conditions to establish the LDP for the DTHP via G\"artner-Ellis Theorem. Therefore, one natural possible direction this work can be extended is to explore the existence of the LDP for the DTHP. Further, one can also study the LDP for the compensator of the DTHP.

\bibliographystyle{elsarticle-harv} 
\bibliography{cas-refs}

\begin{thebibliography}{21}
\expandafter\ifx\csname natexlab\endcsname\relax\def\natexlab#1{#1}\fi
\providecommand{\url}[1]{\texttt{#1}}
\providecommand{\href}[2]{#2}
\providecommand{\path}[1]{#1}
\providecommand{\DOIprefix}{doi:}
\providecommand{\ArXivprefix}{arXiv:}
\providecommand{\URLprefix}{URL: }
\providecommand{\Pubmedprefix}{pmid:}
\providecommand{\doi}[1]{\href{http://dx.doi.org/#1}{\path{#1}}}
\providecommand{\Pubmed}[1]{\href{pmid:#1}{\path{#1}}}
\providecommand{\bibinfo}[2]{#2}
\ifx\xfnm\relax \def\xfnm[#1]{\unskip,\space#1}\fi
\bibitem[{Bacry et~al.(2013)Bacry, Delattre, Hoffmann and Muzy}]{bacry2013some}
\bibinfo{author}{Bacry, E.}, \bibinfo{author}{Delattre, S.}, \bibinfo{author}{Hoffmann, M.}, \bibinfo{author}{Muzy, J.F.}, \bibinfo{year}{2013}.
\newblock \bibinfo{title}{Some limit theorems for hawkes processes and application to financial statistics}.
\newblock \bibinfo{journal}{Stochastic Processes and their Applications} \bibinfo{volume}{123}, \bibinfo{pages}{2475--2499}.
\bibitem[{Bordenave and Torrisi(2007)}]{bordenave2007large}
\bibinfo{author}{Bordenave, C.}, \bibinfo{author}{Torrisi, G.L.}, \bibinfo{year}{2007}.
\newblock \bibinfo{title}{Large deviations of poisson cluster processes}.
\newblock \bibinfo{journal}{Stochastic Models} \bibinfo{volume}{23}, \bibinfo{pages}{593--625}.
\bibitem[{Br{\'e}maud and Massouli{\'e}(1996)}]{bremaud1996stability}
\bibinfo{author}{Br{\'e}maud, P.}, \bibinfo{author}{Massouli{\'e}, L.}, \bibinfo{year}{1996}.
\newblock \bibinfo{title}{Stability of nonlinear hawkes processes}.
\newblock \bibinfo{journal}{The Annals of Probability} , \bibinfo{pages}{1563--1588}.
\bibitem[{Daley et~al.(2003)Daley, Vere-Jones et~al.}]{daley2003introduction}
\bibinfo{author}{Daley, D.J.}, \bibinfo{author}{Vere-Jones, D.}, et~al., \bibinfo{year}{2003}.
\newblock \bibinfo{title}{An introduction to the theory of point processes: volume I: elementary theory and methods}.
\newblock \bibinfo{publisher}{Springer}.
\bibitem[{Dassios and Zhao(2011)}]{dassios2011dynamic}
\bibinfo{author}{Dassios, A.}, \bibinfo{author}{Zhao, H.}, \bibinfo{year}{2011}.
\newblock \bibinfo{title}{A dynamic contagion process}.
\newblock \bibinfo{journal}{Advances in applied probability} \bibinfo{volume}{43}, \bibinfo{pages}{814--846}.
\bibitem[{Errais et~al.(2010)Errais, Giesecke and Goldberg}]{errais2010affine}
\bibinfo{author}{Errais, E.}, \bibinfo{author}{Giesecke, K.}, \bibinfo{author}{Goldberg, L.R.}, \bibinfo{year}{2010}.
\newblock \bibinfo{title}{Affine point processes and portfolio credit risk}.
\newblock \bibinfo{journal}{SIAM Journal on Financial Mathematics} \bibinfo{volume}{1}, \bibinfo{pages}{642--665}.
\bibitem[{Guo and Zeng(2008)}]{guo2008intensity}
\bibinfo{author}{Guo, X.}, \bibinfo{author}{Zeng, Y.}, \bibinfo{year}{2008}.
\newblock \bibinfo{title}{{Intensity process and compensator: A new filtration expansion approach and the Jeulin–Yor theorem}}.
\newblock \bibinfo{journal}{The Annals of Applied Probability} \bibinfo{volume}{18}, \bibinfo{pages}{120 -- 142}.
\bibitem[{Hawkes(1971)}]{hawkes1971spectra}
\bibinfo{author}{Hawkes, A.G.}, \bibinfo{year}{1971}.
\newblock \bibinfo{title}{Spectra of some self-exciting and mutually exciting point processes}.
\newblock \bibinfo{journal}{Biometrika} \bibinfo{volume}{58}, \bibinfo{pages}{83--90}.
\bibitem[{Ivanoff et~al.(2007)Ivanoff, Merzbach and Plante}]{ivanoff2007compensator}
\bibinfo{author}{Ivanoff, B.G.}, \bibinfo{author}{Merzbach, E.}, \bibinfo{author}{Plante, M.}, \bibinfo{year}{2007}.
\newblock \bibinfo{title}{A compensator characterization of point processes on topological lattices}.
\newblock \bibinfo{journal}{Electron. J. Probab.} \bibinfo{volume}{12}, \bibinfo{pages}{no. 2, 47--74}.
\bibitem[{Selvamuthu et~al.(2023)Selvamuthu, Pandey and Tardelli}]{selvamuthu2023limit}
\bibinfo{author}{Selvamuthu, D.}, \bibinfo{author}{Pandey, S.}, \bibinfo{author}{Tardelli, P.}, \bibinfo{year}{2023}.
\newblock \bibinfo{title}{Limit theorems for an extended inverse hawkes process with general exciting functions}.
\newblock \bibinfo{journal}{Statistics \& Probability Letters} \bibinfo{volume}{197}, \bibinfo{pages}{109817}.
\bibitem[{Seol(2015)}]{seol2015limit}
\bibinfo{author}{Seol, Y.}, \bibinfo{year}{2015}.
\newblock \bibinfo{title}{Limit theorems for discrete hawkes processes}.
\newblock \bibinfo{journal}{Statistics \& Probability Letters} \bibinfo{volume}{99}, \bibinfo{pages}{223--229}.
\bibitem[{Seol(2017)}]{SEOL2017165}
\bibinfo{author}{Seol, Y.}, \bibinfo{year}{2017}.
\newblock \bibinfo{title}{Limit theorems for the compensator of hawkes processes}.
\newblock \bibinfo{journal}{Statistics \& Probability Letters} \bibinfo{volume}{127}, \bibinfo{pages}{165--172}.
\bibitem[{Seol(2019)}]{seol2019limit}
\bibinfo{author}{Seol, Y.}, \bibinfo{year}{2019}.
\newblock \bibinfo{title}{Limit theorems for an inverse markovian hawkes process}.
\newblock \bibinfo{journal}{Statistics \& Probability Letters} \bibinfo{volume}{155}, \bibinfo{pages}{108580}.
\bibitem[{Wang(2022)}]{wang2022limit}
\bibinfo{author}{Wang, H.}, \bibinfo{year}{2022}.
\newblock \bibinfo{title}{Limit theorems for a discrete-time marked hawkes process}.
\newblock \bibinfo{journal}{Statistics \& Probability Letters} \bibinfo{volume}{184}, \bibinfo{pages}{109368}.
\bibitem[{Wang(2023)}]{wang2023large}
\bibinfo{author}{Wang, H.}, \bibinfo{year}{2023}.
\newblock \bibinfo{title}{Large and moderate deviations for a discrete-time marked hawkes process}.
\newblock \bibinfo{journal}{Communications in Statistics-Theory and Methods} \bibinfo{volume}{52}, \bibinfo{pages}{6037--6062}.
\bibitem[{Xu et~al.(2021)Xu, Zhu and Wang}]{xu2021self}
\bibinfo{author}{Xu, Y.}, \bibinfo{author}{Zhu, L.}, \bibinfo{author}{Wang, H.}, \bibinfo{year}{2021}.
\newblock \bibinfo{title}{A self-and mutual-exciting model for discrete-time data: Case study on online money market fund}.
\newblock \bibinfo{journal}{Available at SSRN 3665436} .
\bibitem[{Zhu(2013a)}]{zhu2013central}
\bibinfo{author}{Zhu, L.}, \bibinfo{year}{2013}a.
\newblock \bibinfo{title}{Central limit theorem for nonlinear hawkes processes}.
\newblock \bibinfo{journal}{Journal of Applied Probability} \bibinfo{volume}{50}, \bibinfo{pages}{760--771}.
\bibitem[{Zhu(2013b)}]{zhu2013moderate}
\bibinfo{author}{Zhu, L.}, \bibinfo{year}{2013}b.
\newblock \bibinfo{title}{Moderate deviations for hawkes processes}.
\newblock \bibinfo{journal}{Statistics \& Probability Letters} \bibinfo{volume}{83}, \bibinfo{pages}{885--890}.
\bibitem[{Zhu(2013c)}]{zhu2013ruin}
\bibinfo{author}{Zhu, L.}, \bibinfo{year}{2013}c.
\newblock \bibinfo{title}{Ruin probabilities for risk processes with non-stationary arrivals and subexponential claims}.
\newblock \bibinfo{journal}{Insurance: Mathematics and Economics} \bibinfo{volume}{53}, \bibinfo{pages}{544--550}.
\bibitem[{Zhu(2014)}]{zhu2014process}
\bibinfo{author}{Zhu, L.}, \bibinfo{year}{2014}.
\newblock \bibinfo{title}{Process-level large deviations for nonlinear hawkes point processes}, in: \bibinfo{booktitle}{Annales de l'IHP Probabilit{\'e}s et statistiques}, pp. \bibinfo{pages}{845--871}.
\bibitem[{Zhu(2015)}]{zhu2015large}
\bibinfo{author}{Zhu, L.}, \bibinfo{year}{2015}.
\newblock \bibinfo{title}{Large deviations for {M}arkovian nonlinear {H}awkes processes}.
\newblock \bibinfo{journal}{Ann. Appl. Probab.} \bibinfo{volume}{25}, \bibinfo{pages}{548--581}.

\end{thebibliography}





\end{document}